\documentclass[english]{article}
\usepackage[T1]{fontenc}
\usepackage[latin1]{inputenc}
\usepackage{amsmath}
\usepackage{amssymb}
\usepackage[authoryear]{natbib}
\usepackage{amsfonts}
\usepackage{babel}

\makeatletter

\newtheorem{theorem}{Theorem}

\newtheorem{corollary}[theorem]{Corollary}

\newtheorem{definition}[theorem]{Definition}

\newtheorem{lemma}[theorem]{Lemma}

\newtheorem{proposition}[theorem]{Proposition}

\makeatother
\input{tcilatex}
\begin{document}

\title{A Limit Theorem for Products of Free Unitary Operators}
\author{Vladislav Kargin \thanks{%
Courant Institute of Mathematical Sciences, 251 Mercer Street, New York, NY
10012; kargin@cims.nyu.edu} \thanks{%
The author would like to express his gratitude to Diana Bloom for her help
with editing, and to Professor Raghu Varadhan for useful discussions.}}
\date{}
\maketitle

\begin{center}
\textbf{Abstract}
\end{center}

\begin{quotation}
This paper establishes necessary and sufficient conditions for the products
of freely independent unitary operators to converge in distribution to the
uniform law on the unit circle.
\end{quotation}



\section{Introduction}

Suppose $X$ is a unitary $n$-by-$n$ matrix. Then $X$ has $n$ eigenvalues,
which are all located on the unit circle. If we give each eigenvalue a
weight of $n^{-1},$ then we can think about the distribution of these
eigenvalues as a probability distribution supported on $n$ points of the
unit circle. More generally, if $X$ is a unitary operator in a finite von
Neumann algebra, then we can define a spectral probability distribution of $%
X,$ which is supported on the unit circle (see, e.g., Section 1.1 in %
\citet{hiai_petz00}).

If we have several unitary operators $X_{1}$, \ldots , $X_{n}$, then it is
natural to ask about the spectral distribution of their product. In general,
we cannot determine this distribution without more information about
relations among operators $X_{1},$ \ldots , $X_{n}.$ However, if $X_{1},$
\ldots , $X_{n}$ are infinite-dimensional and, in a certain sense, in a
general position relative to each other, then the spectral distribution of
their product is computable. The idea of a general position was formalized
by Voiculescu in his concept of freeness of operators (see %
\citet{voiculescu83},\citet{voiculescu86}, and a textbook by %
\citet{hiai_petz00}). If operators $X_{1},$ \ldots , $X_{n}$ are free and
unitary and their spectral probability distributions are $\mu _{1},$ \ldots
, $\mu _{n},$ respectively, then the distribution of their product is
determined uniquely. This distribution is called the free multiplicative
convolution of measures $\mu _{1},$ $\ldots ,$ $\mu _{n}$ and denoted as $%
\mu _{1}\boxtimes \ldots \boxtimes \mu _{n}.$

What can we say about the asymptotic behavior of $\mu ^{\left( n\right) }=:$ 
$\mu _{1}\boxtimes \ldots \boxtimes \mu _{n},$ as $n$ increases to infinity?
In particular, what are necessary and sufficient conditions on $\mu _{i}$
that ensure that $\mu ^{\left( n\right) }$ converges to the uniform
distribution on the unit circle?

To answer this question, let us define the expectation with respect to the
measure $\mu_{i}$. This is a functional that maps functions analytic in a
neighborhood of the unit circle to complex numbers: 
\begin{equation*}
E_{\mu_{i}}f=:\int_{\left\vert \xi\right\vert
=1}f\left(\xi\right)d\mu_{i}\left(\xi\right).
\end{equation*}
If unitary operator $X_{i}$ has the spectral probability distribution $%
\mu_{i}$, then we will also write:%
\begin{equation*}
Ef\left(X_{i}\right)=:E_{\mu_{i}}f.
\end{equation*}
In particular, $EX_{i}$ denotes $\int_{\left\vert \xi\right\vert =1}\xi
d\mu_{i}\left(\xi\right).$ Then the answer is given by the following theorem:

\begin{theorem}
\label{theorem_main_result}Suppose $\left\{ X_{i}\right\} _{i=1}^{\infty }$
are free unitary operators with spectral measures $\mu _{i}.$ The measures $%
\mu ^{\left( n\right) }$ of the products $\Pi _{n}=:X_{1}\ldots X_{n}$
converge to the uniform measure on the unit circle if and only if at least
one of the following situations holds: \newline
(i) There exist two indices $i\neq j$ such that $EX_{i}=EX_{j}=0;$\newline
(ii) There exists exactly one index $i$ such that $EX_{i}=0,$ and $%
\prod\nolimits_{k=i+1}^{n}EX_{k}\rightarrow 0$ as $n\rightarrow \infty ;$%
\newline
(iii) There exists exactly one index $i$ such that $X_{i}$ has the uniform
distribution;\newline
(iv) $EX_{k}\neq 0$ for all $k,$ and $\prod\nolimits_{k=1}^{n}EX_{k}%
\rightarrow 0$ as $n\rightarrow \infty .$
\end{theorem}

In other words, convergence of $\mu^{\left(n\right)}$ to the uniform law
implies that $\prod\nolimits _{k=1}^{n}EX_{k}\rightarrow0,$ and the only
case when the reverse implication fails is when $EX_{i}=0$ for exactly one $%
X_{i},$ the measure $\mu_{i}$ is not uniform, and $\prod\nolimits
_{k=i+1}^{n}EX_{k}\nrightarrow0$ as $n\rightarrow\infty.$ Note that cases
(ii) and (iii) above are not exclusive. It may happen that both $\mu_{i}$ is
uniform and $\prod\nolimits _{k=i+1}^{n}EX_{k}\rightarrow0$ as $%
n\rightarrow\infty.$ In this case, both (ii) and (iii) hold, and $%
\mu^{\left(n\right)}$ converges to the uniform law.

This theorem can be thought of as a limit theorem about free multiplicative
convolutions of measures on the unit circle. There is some literature about
traditional multiplicative convolutions of measures on the unit circle, or
more generally, about convolutions of measures on compact groups. For the
unit circle, this investigation was started by \citet{levy39}. Then it was
continued by \citet{kawada_ito40}, who studied compact groups, and %
\citet{dvoretzky_wolfowitz51} and \citet{vorobev54}, who both considered the
case of commutative finite groups. These researchers found an important
necessary condition for convergence of convolutions to the uniform law. This
condition requires that there should be no normal subgroup such that the
convolved measures are supported entirely in an equivalence class relative
to this subgroup. This condition is sufficient if summands are identically
distributed. If they are not, then there are some sufficient and necessary
conditions, which are especially useful if the group is cyclic. A textbook
presentation with further references can be found in \citet{grenander63}.

Recent investigations of convolutions on groups are mostly concerned with
the speed of convergence of convolved measures to the uniform law. For a
description of progress in this direction, the reader can consult surveys in %
\citet{diaconis88} and \citet{saloff-coste04}.

It turns out that free convolutions converge to the uniform law under much
weaker conditions than usual convolutions. As an example, consider the
distributions that are concentrated on $-1$ and $+1$. Let measure $\mu_{k}$
put the weight $p_{k}\,$\ on $+1.$ Then usual convolutions remain
concentrated on $-1$ and $+1,$ and therefore they have no chance to converge
to the uniform distribution on the unit circle. In contrast, we will show
that free convolutions do converge to the uniform law, provided that either $%
\prod\nolimits _{k=k_{0}}^{n}\left(2p_{k}-1\right)\rightarrow0$ for
arbitrarily large $k_{0},$ or there exist two indices $i$ and $j$ such that $%
p_{i}=p_{j}=1/2.$

The rest of the paper is organized as follows. Section \ref%
{section_background}\ provides the necessary background. In Section \ref%
{section_outline_of_proof} we outline the proof. Section \ref%
{section_auxiliary_results}\ \ derives some auxiliary results that will be
used in the proof. Section \ref{section_analysis} proves the main result
(Theorem \ref{theorem_main_result}). Section \ref{section_key_estimate}
derives the key estimate used in the proof. And Section \ref%
{section_conclusion} concludes.

\section{Definitions and Background\label{section_background}}

\begin{definition}
A \emph{non-commutative probability space} is a pair $\left(\mathcal{A}%
,E\right),$ where $\mathcal{A}$ is a unital $C^{\ast}$-algebra of bounded
linear operators acting on a complex separable Hilbert space and $E$ is a
linear functional from $\mathcal{A}$ to complex numbers. The operators from
algebra $\mathcal{A}$ are called \emph{non-commutative random variables,} or
simply \emph{random variables}, and the functional $E$ is called the \emph{%
expectation.}
\end{definition}

The linear functional $E$ is assumed to satisfy the following properties (in
addition to linearity): i) $E(I)=1;$ ii) $E(A^{\ast})=E(A);$ iii) $%
E(AA^{\ast})\geq0;$ iv) $E(AA^{\ast})=0$ implies $A=0;$ and v) if $%
A_{n}\rightarrow A$ then $E\left(A_{n}\right)\rightarrow E\left(A\right)$,
where convergence of operators is in norm.

If $P\left( d\lambda \right) $ is the spectral resolution associated with a
unitary operator $A,$ then we can define a measure $\mu \left( d\lambda
\right) =E\left( P\left( d\lambda \right) \right) .$ It is easy to check
that $\mu $ is a probability measure supported on the unit circle. We call
this measure, $\mu ,$ the \emph{spectral probability measure associated with
operator} $A$ \emph{and expectation} $E.$

The most important concept in free probability theory is that of free
independence of non-commuting random variables. Let a set of r.v. $%
A_{1,}\ldots ,A_{n}$ be given. With each of them we can associate an algebra 
$\mathcal{A}_{i},$ which is generated by $A_{i}$; that is, it is the closure
of all polynomials in variables $A_{i}$ and $A_{i}^{\ast }.$ Let $\overline{A%
}_{i}$ denote an arbitrary element of algebra $\mathcal{A}_{i}.$

\begin{definition}
The algebras $\mathcal{A}_{1,}\ldots ,\mathcal{A}_{n}$ (and variables $%
A_{1,}\ldots ,A_{n}$ that generate them) are said to be \emph{freely
independent} or \emph{free}, if the following condition holds: 
\begin{equation*}
\varphi \left( \overline{A}_{i(1)}...\overline{A}_{i(m)}\right) =0,
\end{equation*}%
provided that $\varphi \left( \overline{A}_{i(s)}\right) =0$ and $i(s+1)\neq
i(s)$.
\end{definition}

For more information about non-commutative probability spaces and free
operators we refer the reader to Sections 2.2 - 2.5 in the book by %
\citet{voiculescu_dykema_nica92}.

We will use two results regarding the free operators, which we cite without
proofs. The first one is formula (2.2.3) on page 44 in \citet{hiai_petz00}.

\begin{proposition}
\label{joint_moment_reduction}Let $\mathcal{A}_{1},\ldots ,\mathcal{A}_{m}$
be free sub-algebras of $\mathcal{A}$, and let $A_{1},\ldots ,A_{n}$ be a
sequence of random variables, $A_{k}\in \mathcal{A}_{i(k)}$ such that $%
i(k)\neq i(k+1).$ Then 
\begin{equation}
E\left( A_{1}...A_{n}\right) =\sum_{r=1}^{n}\sum_{1\leq k_{1}<...<k_{r}\leq
n}\left( -1\right) ^{r-1}E\left( A_{k_{1}}\right) \ldots E\left(
A_{k_{r}}\right) E\left( A_{1}\ldots \widehat{A}_{k_{1}}\ldots \widehat{A}%
_{k_{r}}\ldots A_{n}\right) ,  \label{free_relation_property1}
\end{equation}%
where $\symbol{94}$ denotes terms that are omitted.\newline
\end{proposition}

\textbf{Remark:} Note that on the right-hand side the expectations are taken
of the products that have no more than $n-1$ terms. So a recursive
application of this formula reduces computation of $E\left( X_{1}\ldots
X_{n}\right) $ to a polynomial in the moments of the individual variables.

The second result is the Voiculescu multiplication theorem. To formulate it
we need some additional definitions.

Define the $\psi$\emph{-function} of a bounded random variable $X$ as 
\begin{equation}
\psi_{X}\left(z\right)=:\sum_{k=1}^{\infty}E\left(X^{k}\right)z^{k}.
\label{definition_psi_function}
\end{equation}
If $X$ is unitary operator with the spectral measure $\mu$, then we can
write:

\begin{equation*}
\psi_{\mu}\left(z\right)=\int_{\left\vert \xi\right\vert =1}\frac{1}{1-\xi z}%
d\mu\left(\xi\right)-1.
\end{equation*}

It is very useful to note that the $\psi$-function is related to the Poisson
transform of the measure $\mu.$ Indeed, since 
\begin{equation*}
\mathrm{Re}\frac{1}{1-\xi z}=\frac{1-r\cos\left(\omega-\theta\right)}{%
1-2r\cos\left(\omega-\theta\right)+r^{2}},
\end{equation*}
where $\xi=e^{-i\theta}$ and $z=re^{i\omega},$ therefore, 
\begin{equation*}
\mathrm{Re}\frac{1}{1-\xi z}=\frac{1}{2}+\pi P\left(r,\omega-\theta\right)
\end{equation*}
where $P\left(r,\theta\right)$ is the \emph{Poisson kernel}:%
\begin{equation*}
P\left(r,\theta\right)=\frac{1}{2\pi}\frac{1-r^{2}}{1-2r\cos\theta+r^{2}}.
\end{equation*}

Recall that the \emph{Poisson transform} of a measure $\mu $ supported on
the unit circle is defined as 
\begin{equation*}
U_{\mu }\left( z\right) =:\int_{-\pi }^{\pi }P\left( r,\omega -\theta
\right) d\mu \left( \theta \right) ,
\end{equation*}%
where $z=re^{i\omega }$. (Here we have identified measures on the unit
circle and on the interval $\left[ -\pi ,\pi \right) $: $\mu \left( d\theta
\right) =\mu \left\{ \xi :\left\vert \xi \right\vert =1\text{ and }\arg \xi
\in d\theta \right\} $). Hence, 
\begin{equation}
U_{\mu }\left( z\right) =\frac{1}{\pi }\mathrm{Re}\text{ }\psi \left(
z\right) +\frac{1}{2\pi }.  \label{Poisson_transform_and_psi}
\end{equation}

Let $\psi_{X}^{-1}\left(u\right)$ denote the functional inverse of $%
\psi_{X}\left(z\right)$ in a neighborhood of $z=0,$ where $%
\psi_{X}\left(z\right)$ is as defined in (\ref{definition_psi_function}).
(This inversion is possible provided that $E\left(X\right)\neq0$.) Define
also 
\begin{equation*}
S_{X}\left(u\right)=\frac{u+1}{u}\psi_{X}^{-1}\left(u\right).
\end{equation*}

\begin{theorem}
{[}Voiculescu] \label{Voiculescu_multiplication}Suppose $X$ and $Y$ are
bounded free random variables. Suppose also that $E\left(X\right)\neq0$ and $%
E\left(Y\right)\neq0.$ Then 
\begin{equation*}
S_{XY}\left(z\right)=S_{X}\left(z\right)S_{Y}\left(z\right).
\end{equation*}
\end{theorem}

The original proof can be found in \citet{voiculescu87}. A simpler proof was
given by \citet{haagerup97}. Using this theorem, it is possible to compute
the free convolution of two measures, $\mu _{1}$ and $\mu _{2}.$ First, we
can compute their $\psi $-functions, $\psi _{\mu _{1}}$ and $\psi _{\mu
_{2}} $. Then we invert them and obtain the $S$-functions, $S_{\mu _{1}}$
and $S_{\mu _{2}}$. Their product is the $S$-function of the free
convolution, $S_{\mu _{1}\boxtimes \mu _{2}},$ and we can compute $\psi
_{\mu _{1}\boxtimes \mu _{2}}$ by inversion. This determines the Poisson
transform of $\mu _{1}\boxtimes \mu _{2},$ from which we can determine the
measure itself. (For the one-to-one relation of Poisson transforms and
corresponding measures, see Theorem I.3.1 on page 15 and a comment on page
20 in \citet{Garnett81}.)

\section{Outline of the Proof\label{section_outline_of_proof}}

Let $\Pi _{n}$ denote the partial products: $\Pi _{n}=X_{1}\ldots X_{n}.$ We
denote $E\left( X_{i}\right) $ as $a_{i},$ and $E\left( \Pi _{n}\right) $ as 
$a_{\left( n\right) }.$ First, note that it is enough to consider the case
when all $a_{i}$ are real and non-negative. Indeed, for an arbitrary
sequence of real constants $\theta _{n},$ the sequence of operators $%
e^{i\theta _{n}}\Pi _{n}$ converges in distribution to the uniform law if
and only if the sequence $\Pi _{n}$ converges in distribution to the uniform
law. (Indeed, if, say, $\Pi _{n}$ does not converge in distribution to the
uniform law, then we can find an integer $k$ such that $\left\vert
\int_{0}^{2\pi }e^{ik\theta }\mu _{\left( n\right) }\left( d\theta \right)
\right\vert \nrightarrow 0,$ where $\mu _{\left( n\right) }$ denotes the
measure of $\Pi _{n}.\,$But then $\left\vert \int_{0}^{2\pi }e^{ik\theta
}e^{i\theta _{n}}\mu _{\left( n\right) }\left( d\theta \right) \right\vert
=\left\vert \int_{0}^{2\pi }e^{ik\theta }\mu _{\left( n\right) }\left(
d\theta \right) \right\vert \nrightarrow 0,$ and this implies $e^{i\theta
_{n}}\Pi _{n}$ does not converge in distribution to the uniform law.)
Therefore if $a_{i}=E\left( X_{i}\right) $ is not real and positive, then we
can replace $X_{i}$ with $e^{-i\arg a_{i}}X_{i}$ without affecting the
convergence of $\Pi _{n}.$

We divide the analysis into the following cases:

Case I $\; a_{\left(n\right)}\nrightarrow0.$

Case II $\; a_{\left(n\right)}\rightarrow0,$ and there are at least two
indices, $i$ and $j,$ such that $a_{i}=a_{j}=0.$

Case III $\; a_{\left(n\right)}\rightarrow0,$ and for all $i,$ $a_{i}>0.$

\qquad Subcase III.1 $\;\lim\inf a_{i}=0.$

\qquad Subcase III.2 $\;\lim\inf a_{i}=\underline{a}>0.$

Case IV $a_{\left(n\right)}\rightarrow0,$ and there exists exactly one index 
$i,$ such that $a_{i}=0.$

We will show that without loss of generality we can assume in this case that 
$a_{1}=0,$ and $a_{k}>0$ for all $k>1.$

\qquad Subcase IV.1 $\; X_{1}$ has the uniform distribution.

\qquad Subcase IV.2 $\; X_{1}$ does not have the uniform distribution, and $%
\prod\nolimits _{k=2}^{n}a_{n}\rightarrow0$ as $n\rightarrow\infty.$

\qquad Subcase IV.3 $\; X_{1}$ does not have the uniform distribution, and $%
\prod\nolimits _{k=2}^{n}a_{n}\nrightarrow0$ as $n\rightarrow\infty.$

We will show that $\Pi_{n}$ does not converge to the uniform law if and only
if either Case I or Case IV.3 holds.

\section{Auxiliary Lemmas \label{section_auxiliary_results}}

We will need to perform functional inversions. A useful tool for doing this
is Lagrange's formula.

\begin{lemma}
{[}Lagrange's inversion formula] \label{Lagranges_series_around0} Suppose
that (i) $f$ is a function of a complex variable $z,$ which is analytic in a
neighborhood of $z=0,$ (ii) $f\left( 0\right) =0,$ and (iii) $f^{\prime
}\left( 0\right) =a\neq 0$. Then the functional inverse of $f\left( z\right) 
$ is well defined in a neighborhood of $0,$ and the Taylor series of the
inverse is given by the following formula:%
\begin{equation*}
f^{-1}\left( u\right) =\frac{u}{a}+\sum_{k=2}^{\infty }\left[ \frac{1}{k}%
\mathrm{res}_{z=0}\frac{1}{f(z)^{k}}\right] u^{k},
\end{equation*}%
where $\mathrm{res}_{z=0}$ denotes the Cauchy residue at $0$. In addition,%
\begin{equation*}
f^{-1}\left( u\right) =\frac{u}{a}+\sum_{k=2}^{\infty }\left[ \frac{1}{2\pi
ik}\oint_{\gamma }\frac{dz}{f(z)^{k}}\right] u^{k},
\end{equation*}%
where $\gamma $ is a circle around $0,$ inside which $f$ has only one zero.
\end{lemma}

For a proof see Section 7.32 in \citet{whittaker_watson27}.

We will also use the lemmas below:

\begin{lemma}
\label{lemma_moment_product_estimate}Suppose $A$ and $B$ are free unitary
operators, $\left\vert E\left( A\right) \right\vert \leq a$ and $\left\vert
E\left( B\right) \right\vert \leq b.$ Then for all integer $k\geq 1,$%
\begin{equation*}
\left\vert E\left[ \left( AB\right) ^{k}\right] \right\vert \leq M_{k}\max
\left( a,b\right)
\end{equation*}%
for certain constants $M_{k},$ which depend only on $k$.
\end{lemma}

\textbf{Proof}: If we expand $E\left[ \left( AB\right) ^{k}\right] $ using
Proposition \ref{joint_moment_reduction}, then we can observe that each term
in the expansion contains either $E\left( A\right) $ or $E\left( B\right) $
as a separate multiple. The remaining multiples in this term are $\leq 1$ in
absolute value; therefore, we can bound each term by $\max \left( a,b\right)
.$ The number of terms in this expansion is bounded by a constant, $M_{k}.$
Therefore, $\left\vert E\left[ \left( AB\right) ^{k}\right] \right\vert $ is
bounded by $M_{k}\max \left( a,b\right) .$ QED.

In the following lemmas we use the fact that the sequence of probability
measures $\mu _{i},$ supported on the unit circle, converges to the uniform
law if and only if all their moments converge to $0,$ that is, iff for each $%
k\geq 1$, $\int_{\left\vert \xi \right\vert =1}\xi ^{k}d\mu _{i}\left( \xi
\right) \rightarrow 0$ as $i\rightarrow \infty .$ For completeness we give a
proof of this result.

Let us define $c_{k}^{\left(i\right)}=:$ $\int_{\left\vert \xi\right\vert
=1}\xi^{k}d\mu_{i}\left(\xi\right).$ Note that for a fixed $i,$ $%
c_{k}^{\left(i\right)}$ are coefficients in the Taylor series of $%
\psi_{i}\left(z\right)$, i.e., the $\psi$-function of the measure $\mu_{i}.$

\begin{lemma}
\label{lemma_psi_convergence_on_compact_subsets}Let $\mu _{i}$ be a sequence
of measures supported on the unit circle. If for each $k$ the coefficients $%
c_{k}^{\left( i\right) }\rightarrow 0$ as $i\rightarrow \infty ,$ then $\psi
_{i}\left( z\right) \rightarrow 0$ uniformly on compact subsets of the open
unit disc.
\end{lemma}

\textbf{Proof}: Let $\Omega $ be a compact subset of the open unit disc, and
let $\Omega \subset D_{r},$ where $D_{r}$ denotes a closed disc with the
radius $r<1.$ Fix an $\varepsilon \in \left( 0,1\right) .$ Then we can find
such a $k_{0}$ that 
\begin{equation*}
\left\vert \sum_{k=k_{0}}^{\infty }c_{k}^{\left( j\right) }z^{k}\right\vert
<\varepsilon /2
\end{equation*}%
for all $z\in D_{r}$ and all $j.$ Indeed, $\left\vert c_{k}^{\left( j\right)
}\right\vert \leq 1,$ and therefore, 
\begin{equation*}
\left\vert \sum_{k=k_{0}}^{\infty }c_{k}^{\left( j\right) }z^{k}\right\vert
\leq \frac{r^{k_{0}}}{1-r},
\end{equation*}%
so we can take $k_{0}$ to be any integer greater than or equal to $\log
(\varepsilon \left( 1-r\right) /2)/\log r$.

Given $k_{0},$ we choose a $j_{0}$ so large that for all $j>j_{0}$ and all $%
k<k_{0},$ we have $\left\vert c_{k}^{\left( j\right) }\right\vert
<\varepsilon /\left( 2k_{0}\right) .$ This is possible because by assumption
for each $k$ coefficients $c_{k}^{\left( j\right) }$ converge to zero as $%
j\rightarrow \infty ,$ and we consider only a fixed finite number of
possible $k.$

Consequently,%
\begin{equation*}
\left\vert \sum_{k=1}^{k_{0}-1}c_{k}^{\left(j\right)}z^{k}\right\vert
\leq\sum_{k=1}^{k_{0}-1}\left\vert c_{k}^{\left(j\right)}\right\vert
<\varepsilon/2
\end{equation*}
for every $j>j_{0}$ and all $z\in D_{r}.$ Therefore, 
\begin{equation*}
\left\vert \sum_{k=1}^{\infty}c_{k}^{\left(j\right)}z^{k}\right\vert
<\varepsilon
\end{equation*}
for every $j>j_{0}$ and all $z\in D_{r}.$ Therefore, $\psi_{j}\left(z\right)%
\rightarrow0$ uniformly on $D_{r},$ and therefore on $\Omega.$ Since $\Omega$
was arbitrary, we have proved that $\psi_{j}\left(z\right)\rightarrow0$
uniformly on compact subsets of the unit disc. QED.

The fact that $\psi_{j}\left(z\right)\rightarrow0$ implies that the Poisson
transforms of measures $\mu_{j}$ converge to $\frac{1}{2\pi},$ and therefore 
$\mu_{j}\rightarrow\nu$, where $\nu$ is the uniform measure on the unit
disc. Indeed, we only need to invoke the following result:

\begin{proposition}
\label{proposition_Poisson_convergence_on_unit_disc}If Poisson transforms $%
U_{\mu_{j}}\left(z\right)\rightarrow1/\left(2\pi\right)$ uniformly on
compact subsets of the unit disc, then $\mu_{j}$ weakly converges to $\nu,$
where $\nu$ is the uniform probability measure on the unit circle.
\end{proposition}

\textbf{Proof:} This proposition directly follows from Theorem I.3.1 on page
15 in \citet{Garnett81}, adapted to the case of measures on the unit disc.
QED.

\begin{lemma}
\label{lemma_continuity_convergence}Suppose $\left\{ A_{n}\right\}
_{n=1}^{\infty}$ is a sequence of unitary operators that converges in
distribution to the uniform law. Let $\left\{ B_{n}\right\} _{n=1}^{\infty}$
be another sequence of unitary operators, and let the operator $B_{n}$ be
free of$\ $the operator $A_{n}$ for every $n.$ Then the sequence of products 
$B_{n}A_{n}$ converges in distribution to the uniform law. Also, the
sequence $A_{n}B_{n}$ converges to the uniform law.
\end{lemma}

\textbf{Proof:} Let $a_{k}^{\left( n\right) }=:E\left( \left( A_{n}\right)
^{k}\right) .$ By assumption, for each fixed $k,$ the moment $a_{k}^{\left(
n\right) }\rightarrow 0$ as $n\rightarrow \infty .$ If we represent $E\left(
\left( B_{n}A_{n}\right) ^{k}\right) $ as a polynomial in individual moments
of $B_{n}$ and $A_{n}$, then all terms of this polynomial contain at least
one of the moments $a_{i}^{\left( n\right) },$ $i\leq k,$ which are perhaps
multiplied by some other moments. All of these other moments are less than $%
1 $ in absolute value. Therefore, we can write the following estimate: 
\begin{equation*}
\left\vert E\left( \left( B_{n}A_{n}\right) ^{k}\right) \right\vert \leq
M_{k}^{\prime }\max_{i\leq k}\left\{ a_{i}^{\left( n\right) }\right\} ,
\end{equation*}%
where $M_{k}^{\prime }$ is the number of terms in the polynomial. If $k$ is
fixed and $n$ is growing, then the assumption that $A_{n}$ converges in
distribution to the uniform law implies that $\max_{i\leq k}\left\{
a_{i}^{\left( n\right) }\right\} $ converges to zero. Therefore, all moments
of $B_{n}A_{n}$ converge to zero as $n\rightarrow \infty $, and therefore,
by Lemma \ref{lemma_psi_convergence_on_compact_subsets} and Proposition \ref%
{proposition_Poisson_convergence_on_unit_disc}, the sequence $B_{n}A_{n}$
converges in distribution to the uniform law. \ A similar argument proves
that $A_{n}B_{n}$ converges in distribution to the uniform law. QED.

\begin{lemma}
\label{lemma_product_nonconvergence}Suppose that $B$ is a unitary operator, $%
\left\{ A_{n}\right\} $ is a sequence of unitary operators, $B$ is free from
each of $A_{n,}$ $E\left(B\right)\neq0,$ and the sequence $A_{n}$ does not
converge to uniform law. Then the sequence of products $BA_{n}$ does not
converge to the uniform law.
\end{lemma}

\textbf{Proof:} The condition that the sequence $A_{n}$ does not converge to
the uniform law means that for some fixed $k$ the sequence of $k$-th moments
of $A_{n}$ does not converge to zero as $n\rightarrow \infty $. Let $k$ be
the smallest of these indices. By selecting a subsequence we can assume that 
$\left\vert E\left( A_{n}^{k}\right) \right\vert >\alpha >0$ for all $n.$
Consider $E\left( \left( BA_{n}\right) ^{k}\right) $:%
\begin{equation*}
E\left( \left( BA_{n}\right) ^{k}\right) =\left[ E\left( B\right) \right]
^{k}E\left( A_{n}^{k}\right) +\ldots ,
\end{equation*}%
The number of the terms captured by $\ldots $ is finite and depends only on $%
k$.\ Each of these terms includes at least one of $E\left( A_{n}^{i}\right) $
where $i<k,$ and other multipliers in this term are less than 1 in absolute
value. Therefore, each of these terms converges to zero. Hence, for any $%
\varepsilon >0,$ there exist such $N$ that for all $n>N,$ the sum of the
terms captured by $\ldots $ is less than $\varepsilon $ in absolute value.
Take $\varepsilon =\left\vert E\left( B\right) \right\vert ^{k}\alpha /2.$
Then for $n>N,$ we have:%
\begin{equation*}
\left\vert E\left( \left( BA_{n}\right) ^{k}\right) \right\vert \geq
\left\vert E\left( B\right) \right\vert ^{k}\alpha /2.
\end{equation*}%
Therefore, the sequence of products $BA_{n}$ does not converge to the
uniform law. QED.

\begin{lemma}
\label{lemma_product_nonconvergence_2}Suppose that $B$ is a unitary random
variable, $\left\{ A_{n}\right\} $ is a sequence of unitary random
variables, $B$ is free from each of $A_{n,}$ $B$ is not uniform$,$ and the
sequence of expectations $E\left(A_{n}\right)$ does not converge to zero.
Then the sequence of products $BA_{n}$ does not converge to the uniform law.
\end{lemma}

\textbf{Proof:} By selecting a subsequence we can assume that $\left\vert
E\left( A_{n}\right) \right\vert >\alpha >0$ for all $n.$ The assumption
that $B$ is not uniform means that for some $k\geq 1$, $\ E\left(
B^{k}\right) \neq 0.$ Let $k$ be the smallest of such $k.$ Consider $E\left(
\left( BA_{n}\right) ^{k}\right) $:%
\begin{equation*}
E\left( \left( BA_{n}\right) ^{k}\right) =\left[ E\left( A_{n}\right) \right]
^{k}E\left( B^{k}\right) +\ldots ,
\end{equation*}%
Each of the terms in $\ldots $ includes one of $E\left( B^{i}\right) $ where 
$i<k$. Therefore, all terms in $\ldots $ are zero. Hence,%
\begin{equation*}
\left\vert E\left( \left( BA_{n}\right) ^{k}\right) \right\vert =\left\vert 
\left[ E\left( A_{n}\right) \right] ^{k}E\left( B^{k}\right) \right\vert
>\alpha ^{k}\left\vert E\left( B^{k}\right) \right\vert .
\end{equation*}%
Therefore, the sequence of products $BA_{n}$ does not converge to the
uniform law. QED.

\section{Analysis \label{section_analysis}}

We use the following notation: $\psi_{i}$ and $S_{i}$ denote $\psi$- and $S$%
-functions for variables $X_{i}$ (and measures $\mu_{i}$)$,$ and $%
\psi_{\left(n\right)}$ and $S_{\left(n\right)}$ denote these functions for
variables $\Pi_{n}$ (and measures $\mu^{\left(n\right)}$).

\textbf{Case I:} $a_{\left(n\right)}\nrightarrow0$\textbf{. }

Since $E\left(\Pi_{n}\right)=a_{\left(n\right)},$ therefore, if $%
a_{\left(n\right)}\nrightarrow0,$ then $E\left(\Pi_{n}\right)\nrightarrow0.$
Hence, $\Pi_{n}$ cannot converge to the uniform measure on the unit circle.

\textbf{Case II} $a_{\left(n\right)}\rightarrow0,$\textbf{\ and there are at
least two indices} $i$\textbf{\ and} $j$\textbf{\ such that} $a_{i}=a_{j}=0.$

Assume without loss of generality that $j>i.$ Consider $\Pi _{n}$ with $%
n\geq j$ and define $X=:X_{1}\ldots X_{i}$ and $Y=:X_{i+1}\ldots X_{n}$.
Then $\Pi _{n}=XY,$ and $E\left( Y\right) =E\left( X\right) =0.$ Using Lemma %
\ref{lemma_moment_product_estimate}, we obtain that $\left\vert E\left[
\left( \Pi _{n}\right) ^{k}\right] \right\vert =0$ for every $k>0.$
Therefore, the $\psi $-function of $\Pi _{n}$ is zero, and $\Pi _{n}$ has
the uniform distribution on the unit circle.

\textbf{Case III} $a_{\left(n\right)}\rightarrow0,$\textbf{\ and for all} $%
i, $\textbf{\ }$a_{i}>0.$

\textbf{\qquad Subcase III.1} $\lim\inf a_{i}=0.$

In this case we can find a subsequence $a_{n_{i}}$ that monotonically
converges to zero.

Now, consider $\Pi _{j}$, where $j\in \left[ n_{i},n_{i+1}\right) .$ Then we
can write $\Pi _{j}=XY,$ where $X=X_{1}\ldots X_{n_{i}-1},$ and $%
Y=X_{n_{i}}\ldots X_{j}.$ Then $EX\leq a_{n_{i-1}}$ and $EY\leq
a_{n_{i}}\leq a_{n_{i-1}}.$

Applying Lemma \ref{lemma_moment_product_estimate} we get 
\begin{equation*}
\left\vert E\left(\Pi_{j}^{k}\right)\right\vert \leq M_{k}a_{n_{i-1}}.
\end{equation*}
This implies that for a fixed $k,$ $\left\vert
E\left(\Pi_{j}^{k}\right)\right\vert $ approaches zero as $%
j\rightarrow\infty.$ By Lemma \ref{lemma_psi_convergence_on_compact_subsets}
and Proposition \ref{proposition_Poisson_convergence_on_unit_disc}, this
establishes that $\Pi_{j}$ converges to the uniform law.

\textbf{Case III} $a_{\left(n\right)}\rightarrow0,$\textbf{\ and for all} $%
i, $\textbf{\ }$a_{i}>0$\textbf{\ }

\textbf{\qquad Subcase III.2} $\lim\inf a_{i}=\underline{a}>0$.

Let us choose such an $a$ that $0<a<\underline{a}.$ Starting from some $%
j_{0},$ $a_{j}\in \left( a,1\right) .$ Let $\widetilde{\Pi }%
_{n}=X_{j_{0}}\ldots X_{n+j_{0}-1}.$ Then, by Lemmas \ref%
{lemma_continuity_convergence} and \ref{lemma_product_nonconvergence}, $%
\widetilde{\Pi }_{n}$ converges to the uniform law if and only if $\Pi _{n}$
converges to the uniform law \ Hence, without loss of generality we can
restrict our attention to the case when $a_{k}\in \left( a,1\right) $ for
all $k.$

\begin{lemma}
\label{lemma_products_and_sums} Suppose $1\geq a_{k}>0$ for all $k,$ and let 
$\alpha_{i}=:1-a_{i}$. Then $\prod\nolimits _{i=1}^{n}a_{i}\rightarrow0$ if
and only if $\sum_{i=1}^{n}\alpha_{i}\rightarrow\infty.$
\end{lemma}

This is a standard result. For a proof see Section 2.7 in %
\citet{whittaker_watson27}.

Since $\log\left(1-\alpha_{i}\right)\leq-\alpha_{i},$ we also have the
following estimate, which we will find useful later.%
\begin{equation}
\prod\limits _{i=1}^{n}a_{i}\leq\exp\left(-\sum_{i=1}^{n}\alpha_{i}\right).
\label{inequality_product_sum}
\end{equation}

To prove convergence to the uniform law, we have to establish that for every 
$k>0$ the coefficient $c_{k}^{\left( n\right) }$ in the Taylor expansion of
function $\psi _{\left( n\right) }\left( z\right) $ approaches zero as $%
n\rightarrow \infty .$ We know from Lemma \ref{Lagranges_series_around0}
that 
\begin{equation*}
kc_{k}^{\left( n\right) }=\mathrm{res}_{z=0}\frac{1}{\left[ \psi _{\left(
n\right) }^{-1}\left( z\right) \right] ^{k}};
\end{equation*}%
therefore, our main task is to estimate this residual. This is the same as
estimating the coefficient before the term $z^{k-1}$ in the Taylor expansion
of%
\begin{equation*}
f\left( z\right) =\left[ \frac{z}{\psi _{\left( n\right) }^{-1}\left(
z\right) }\right] ^{k}.
\end{equation*}%
We will approach this problem by using the Cauchy inequality (see Section
5.23 in \citet{whittaker_watson27}). Applied to the coefficient before $%
z^{k-1}$ in the Taylor expansion of $f\left( z\right) ,$ this inequality
says that 
\begin{equation}
\left\vert kc_{k}^{\left( n\right) }\right\vert \leq \frac{M\left( r\right) 
}{r^{k-1}},  \label{formula_Cauchy_inequality}
\end{equation}%
where $r>0$ is such that $f\left( z\right) $ is analytic inside $\left\vert
z\right\vert =r,$ and 
\begin{equation*}
M\left( r\right) =:\max_{\left\vert z\right\vert =r}\left\vert f\left(
z\right) \right\vert .
\end{equation*}

It is easy to check that the constant in the Taylor expansion of $z/\psi
_{\left( n\right) }^{-1}\left( z\right) $ is $a_{\left( n\right) }$. So $%
M\left( 0\right) =a_{\left( n\right) },$ which approaches zero as $%
n\rightarrow \infty .$ The main question is how large we can take $r$, so
that $M\left( r\right) $ remains relatively small. In other words, we want
to minimize the right-hand side of (\ref{formula_Cauchy_inequality}) by a
suitable choice of $r.$

\begin{proposition}
\label{proposition_ck_estimate}Suppose that $EX_{i}=a_{i}>a$ for each $i$
and that $a_{\left(n\right)}=:\prod\nolimits _{i=1}^{n}a_{i}\rightarrow0.$
Let $\alpha_{i}=1-a_{i}$. Then for all sufficiently large $n,$ the following
inequality holds:%
\begin{equation*}
\left\vert c_{k}^{\left(n\right)}\right\vert \leq\left(\frac{C}{a^{2}}%
\right)^{k}\left[\left(\sum\limits
_{i=1}^{n}\alpha_{i}\right)\exp\left(-\sum\limits _{i=1}^{n}\alpha_{i}\right)%
\right]^{k},
\end{equation*}
where $C=2^{17}.$
\end{proposition}

\textbf{Proof:} The main tool in the proof is the following proposition:

\begin{proposition}
\label{proposition_f_estimate} Suppose that $\alpha_{i}=:1-a_{i}<1-a$ for
each $i,$ and that $z$ and $n$ are such that 
\begin{equation*}
\left\vert z\right\vert \leq\frac{a^{2}}{6684}\min\left\{
1,\left(\sum_{i=1}^{n}\alpha_{i}\right)^{-1}\right\} .
\end{equation*}
Then,%
\begin{equation*}
\left\vert \frac{z}{\psi_{\left(n\right)}^{-1}\left(z\right)}\right\vert
^{k}\leq\left(2e^{2}\right)^{k}\left(\prod\nolimits
_{i=1}^{n}a_{i}\right)^{k}.
\end{equation*}
\end{proposition}

We will prove this proposition in the next section and assume for now that
it holds.

Let $n_{0}$ be so large that $\sum\nolimits _{i=1}^{n_{0}}\alpha_{i}>1.$ (We
can find such $n_{0}$ because by Lemma \ref{lemma_products_and_sums}, $%
\sum\nolimits _{i=1}^{n}\alpha_{i}\rightarrow\infty$ as $n\rightarrow\infty.$%
) In particular, this implies that $\sum\nolimits _{i=1}^{n}\alpha_{i}>1$
for every $n\geq n_{0}.$ Define $r_{n}=:$ $a^{2}\left(\sum_{i=1}^{n}%
\alpha_{i}\right)^{-1}/6684.$ Then, using Proposition \ref%
{proposition_f_estimate} and formulas (\ref{formula_Cauchy_inequality}) and (%
\ref{inequality_product_sum}), we get: 
\begin{eqnarray*}
\left\vert kc_{k}^{\left(n\right)}\right\vert & \leq &
\left(2e^{2}\right)^{k}\left(\prod\nolimits _{i=1}^{n}a_{i}\right)^{k}\left(%
\frac{6684}{a^{2}}\sum\nolimits _{i=1}^{n}\alpha_{i}\right)^{k-1} \\
& \leq & \left[\frac{2^{17}}{a^{2}}\left(\sum\nolimits
_{i=1}^{n}\alpha_{i}\right)\exp\left(-\sum\nolimits
_{i=1}^{n}\alpha_{i}\right)\right]^{k},
\end{eqnarray*}
provided that $n\geq n_{0}.$ QED.

Using Lemma \ref{lemma_products_and_sums}, we get the following Corollary:

\begin{corollary}
If the assumptions of Proposition \ref{proposition_ck_estimate} hold, then
for each $k,$ the coefficient $c_{k}^{\left(n\right)}\rightarrow0$ as $%
n\rightarrow\infty.$
\end{corollary}

This Corollary shows that in Case III.2 the product $\Pi_{n}$ converges to
the uniform law.

\textbf{Case IV} $a_{\left(n\right)}\rightarrow0,$\textbf{\ and there exists
exactly one index} $i,$\textbf{\ such that} $a_{i}=0.$

First, we want to show that without loss of generality we can assume in this
case that $a_{1}=0,$ and $a_{k}>0$ for all $k>1.$ Indeed, suppose $a_{i}=0$
for $i>1,$ and $a_{j}>0$ for $j<i.$ Let $X=X_{1}\ldots X_{i-1}$ and let $%
\widetilde{\Pi }_{n}=X_{i}\ldots X_{i+n-1}.$ Then $E\left( X\right) \neq 0,$
and using Lemmas \ref{lemma_continuity_convergence} and \ref%
{lemma_product_nonconvergence}, we conclude that $\Pi _{n}$ converges to the
uniform law if and only if $\widetilde{\Pi }_{n}$ converges to the uniform
law.

\textbf{Subcase IV.1} $X_{1}$\textbf{\ has the uniform distribution.}

In this case all moments of $X_{1}$ are zero, i.e., $E\left(
X_{1}^{k}\right) =0$ for all $k>0,$ and Proposition \ref%
{joint_moment_reduction} implies that all moments of $\Pi _{n}$ are zero.
Therefore, $\Pi _{n}$ is uniform for all $n.$

\textbf{Subcase IV.2} $X_{1}$\textbf{\ does not have the uniform
distribution, and} $\prod\nolimits _{k=2}^{n}a_{n}\rightarrow0$\textbf{\ as} 
$n\rightarrow\infty.$

By Case III, the product $X_{2}\ldots X_{n}$ converges to the uniform law,
and using Lemma \ref{lemma_continuity_convergence}, we conclude that $\Pi
_{n}$ also converges to the uniform law.

\textbf{Subcase IV.3} $X_{1}$ \textbf{does not have the uniform distribution
and} $\prod\nolimits _{k=2}^{n}a_{n}\nrightarrow0$\textbf{\ as} $%
n\rightarrow\infty.$

Applying Lemma \ref{lemma_product_nonconvergence_2} to $B=X_{1}$ and $%
A=X_{2}\ldots X_{n},$ we conclude that $\Pi _{n}$ does not converge to the
uniform law.

\section{Proof of Proposition \protect\ref{proposition_f_estimate} \label%
{section_key_estimate}}

Let 
\begin{equation*}
f(z)=:\left(\frac{z}{\psi_{\left(n\right)}^{-1}\left(z\right)}\right)^{k}.
\end{equation*}
Using Theorem \ref{Voiculescu_multiplication}, we can write this function as
follows:%
\begin{equation}
f(z)=\left(\frac{z^{n}}{\left(1+z\right)^{n-1}}\prod_{i=1}^{n}\frac{1}{%
\psi_{i}^{-1}\left(z\right)}\right)^{k}.  \label{formula_product}
\end{equation}

We want to estimate $\left\vert f\left(z\right)\right\vert $ for all
sufficiently small $z.$ We start with some auxiliary estimates, which will
later allow us to estimate $\psi_{i}\left(z\right),$ and then $%
\psi_{i}^{-1}\left(z\right)$ for small $z.$

\begin{lemma}
Suppose $\mu$ is a probability measure on $\left[-\pi,\pi\right)$ such that 
\begin{equation}
\left\vert
\int_{-\pi}^{\pi}\left(e^{i\theta}-1\right)d\mu\left(\theta\right)\right%
\vert \leq\alpha.  \label{condition_expectation}
\end{equation}
Then, i) 
\begin{equation*}
\int_{-\pi}^{\pi}\theta^{2}d\mu\left(\theta\right)\leq\frac{\pi^{2}}{2}%
\alpha<5\alpha;
\end{equation*}
ii) 
\begin{equation*}
\left\vert \int_{-\pi}^{\pi}\theta d\mu\left(\theta\right)\right\vert
\leq\left(1+\frac{\pi^{3}}{12}\right)\alpha<4\alpha,\text{ and}
\end{equation*}
iii) if $k>2,$ then 
\begin{equation*}
\int_{-\pi}^{\pi}\left\vert \theta\right\vert ^{k}d\mu\left(\theta\right)\leq%
\frac{\pi^{k}}{2}\alpha.
\end{equation*}
\end{lemma}

\textbf{Proof:} Condition (\ref{condition_expectation}) implies that 
\begin{equation*}
\int_{-\pi}^{\pi}\left(1-\cos\left(\theta\right)\right)d\mu\left(\theta%
\right)\leq\alpha
\end{equation*}
and that 
\begin{equation*}
\left\vert
\int_{-\pi}^{\pi}\sin\left(\theta\right)d\mu\left(\theta\right)\right\vert
\leq\alpha.
\end{equation*}
Since $1-\cos\theta\geq\left(2/\pi^{2}\right)\theta^{2},$ from the first of
these inequalities we infer that: 
\begin{equation*}
\int_{-\pi}^{\pi}\theta^{2}d\mu\left(\theta\right)\leq\left(\pi^{2}/2\right)%
\alpha,
\end{equation*}
which proves claim i) of the lemma.

Next, note that $\left\vert \sin\theta-\theta\right\vert \leq\left\vert
\theta\right\vert ^{3}/6,$ and that 
\begin{equation*}
\frac{1}{6}\int_{-\pi}^{\pi}\left\vert \theta\right\vert
^{3}d\mu\left(\theta\right)\leq\frac{\pi}{6}\int_{-\pi}^{\pi}\theta^{2}d\mu%
\left(\theta\right)\leq\frac{\pi^{3}}{12}\alpha.
\end{equation*}
Therefore, 
\begin{eqnarray*}
\left\vert \int_{-\pi}^{\pi}\theta d\mu\left(\theta\right)\right\vert & \leq
& \left\vert
\int_{-\pi}^{\pi}\sin\left(\theta\right)d\mu\left(\theta\right)\right\vert
+\left\vert
\int_{-\pi}^{\pi}\left(\theta-\sin\left(\theta\right)\right)d\mu\left(\theta%
\right)\right\vert \\
& \leq & \alpha+\frac{\pi^{3}}{12}\alpha.
\end{eqnarray*}
This proves claim ii) of the lemma.

For claim iii), note that 
\begin{equation*}
\int_{-\pi}^{\pi}\left\vert \theta\right\vert
^{k}d\mu\left(\theta\right)\leq\pi^{k-2}\int_{-\pi}^{\pi}\left\vert
\theta\right\vert ^{2}d\mu\left(\theta\right)\leq\frac{\pi^{k}}{2}\alpha.
\end{equation*}

QED.

\begin{lemma}
\label{lemma_expectations_bound}Suppose Condition (\ref%
{condition_expectation}) holds, and $k$ is a positive integer. Then%
\begin{equation*}
\left\vert
\int_{-\pi}^{\pi}\left(e^{ik\theta}-1\right)d\mu\left(\theta\right)\right%
\vert \leq7k^{3}\alpha.
\end{equation*}
\end{lemma}

\textbf{Proof:} First, remark that $1-\cos\left(k\theta\right)\leq\left(k%
\theta\right)^{2}/2$ and therefore 
\begin{eqnarray*}
\left\vert \int_{-\pi}^{\pi}\left(\cos
k\theta-1\right)d\mu\left(\theta\right)\right\vert & \leq & \frac{k^{2}}{2}%
\int_{-\pi}^{\pi}\theta^{2}d\mu\left(\theta\right) \\
& \leq & \frac{\pi^{2}k^{2}}{4}\alpha.
\end{eqnarray*}

Next, we will use $\left\vert \sin\left(k\theta\right)-k\theta\right\vert
\leq\left(k\left\vert \theta\right\vert \right)^{3}/6$ and write 
\begin{eqnarray*}
\left\vert
\int_{-\pi}^{\pi}\sin\left(k\theta\right)d\mu\left(\theta\right)\right\vert
& \leq & \left\vert \int_{-\pi}^{\pi}k\theta
d\mu\left(\theta\right)\right\vert +\left\vert \frac{1}{6}%
\int_{-\pi}^{\pi}\left(k\left\vert \theta\right\vert
\right)^{3}d\mu\left(\theta\right)\right\vert \\
& \leq & k\left(1+\frac{\pi^{3}}{12}\right)\alpha+\frac{1}{6}k^{3}\frac{%
\pi^{3}}{2}\alpha \\
& \leq & k^{3}\left(1+\frac{\pi^{3}}{6}\right)\alpha.
\end{eqnarray*}

Consequently, 
\begin{eqnarray*}
\left\vert
\int_{-\pi}^{\pi}\left(e^{ik\theta}-1\right)d\mu\left(\theta\right)\right%
\vert & \leq & \alpha\sqrt{\frac{\pi^{4}k^{4}}{16}+k^{6}\left(1+\frac{\pi^{3}%
}{6}\right)^{2}} \\
& \leq & 7k^{3}\alpha.
\end{eqnarray*}
QED.

\begin{lemma}
\label{lemma_psi_approximation}Let $X$ be unitary and $EX=a>0$. If $%
\left\vert z\right\vert \leq1/2$ and $1-a\leq\alpha,$ then 
\begin{equation*}
\left\vert \psi_{X}\left(z\right)-\frac{az}{1-z}\right\vert
\leq716\alpha\left\vert z\right\vert ^{2}.
\end{equation*}
\end{lemma}

\textbf{Proof:} We can write:%
\begin{equation*}
\psi _{X}\left( z\right) -\frac{az}{1-z}=\sum_{k=2}^{\infty }\left( E\left(
X^{k}\right) -a\right) z^{k}.
\end{equation*}%
Therefore, using Lemma \ref{lemma_expectations_bound}, we estimate: 
\begin{eqnarray*}
\left\vert \psi _{X}\left( z\right) -\frac{az}{1-z}\right\vert &\leq
&\sum_{k=2}^{\infty }\left( \left\vert E\left( X^{k}\right) -1\right\vert
+\left\vert 1-a\right\vert \right) z^{k} \\
&\leq &7\alpha \left\vert z\right\vert ^{2}\sum_{k=0}^{\infty }\left[
(k+2)^{3}+1/7\right] \left\vert z\right\vert ^{k} \\
&\leq &716\alpha \left\vert z\right\vert ^{2}.
\end{eqnarray*}%
(Note that 716 is the exact value of the sum $7\sum_{k=0}^{\infty }\left[
(k+2)^{3}+1/7\right] 2^{-k}$.) QED.

To derive a similar estimate for $\psi_{X}^{-1}\left(z\right)$, we need a
couple of preliminary lemmas.

\begin{lemma}
Suppose $X$ is unitary and $EX=a>0.$ Then the function $\psi_{X}\left(z%
\right)$ has only one zero ($z=0$) in the area $\left\vert z\right\vert <a/3$
. If $\left\vert z\right\vert =a/3,$ then $\left\vert
\psi_{X}\left(z\right)\right\vert \geq a^{2}/6.$
\end{lemma}

\textbf{Proof}: Write the following estimate:%
\begin{eqnarray*}
\left\vert \psi _{X}\left( z\right) -az\right\vert &=&\left\vert
\sum_{k=2}^{\infty }E\left( X^{k}\right) z^{k}\right\vert \\
&\leq &\left\vert z\right\vert \sum_{k=1}^{\infty }\left\vert z\right\vert
^{k}=\frac{\left\vert z\right\vert }{1-\left\vert z\right\vert }\left\vert
z\right\vert <\frac{a}{2}\left\vert z\right\vert ,
\end{eqnarray*}%
if $\left\vert z\right\vert <a/3.$ By Rouch$\mathrm{\acute{e}}$'s theorem, $%
\psi _{X}\left( z\right) $ has only one zero in $\left\vert z\right\vert
<a/3.$ The second claim also follows immediately from this estimate. QED.

\begin{lemma}
\label{lemma_psi_X_inverse}Suppose $X$ is unitary and $EX=a>0.$ Then the
function $\psi _{X}^{-1}\left( z\right) $ is analytical for $\left\vert
z\right\vert <a^{2}/6$. If $\left\vert z\right\vert \leq a^{2}/12,$ then 
\begin{equation*}
\left\vert \psi _{X}^{-1}\left( z\right) \right\vert \leq \frac{2}{a}%
\left\vert z\right\vert .
\end{equation*}
\end{lemma}

\textbf{Proof:} Using Lagrange's formula, we can write 
\begin{equation*}
\psi _{X}^{-1}\left( z\right) =\frac{z}{a}+\sum_{k=2}^{\infty }c_{k}z^{k},
\end{equation*}%
where 
\begin{equation*}
c_{k}=\frac{1}{2\pi i}\frac{1}{k}\oint_{\gamma }\frac{du}{\left[ \psi
_{X}\left( u\right) \right] ^{k}}.
\end{equation*}%
By the previous lemma, we can use the circle with the center at $0$ and
radius $a/3$ as $\gamma $, and then we can estimate $c_{k}$ as follows:%
\begin{equation}
\left\vert c_{k}\right\vert \leq \frac{a/3}{k\left( a^{2}/6\right) ^{k}}=%
\frac{2}{ka}\left( \frac{6}{a^{2}}\right) ^{k-1}.
\label{eq:coef_psi_inv_inequality}
\end{equation}%
It follows that the power series for $\psi _{X}^{-1}\left( z\right) $
converges in $\left\vert z\right\vert <a^{2}/6.$ If $\left\vert z\right\vert
<a^{2}/12,$ then we can estimate $\psi _{X}^{-1}\left( z\right) $:%
\begin{eqnarray*}
\left\vert \psi _{X}^{-1}\left( z\right) \right\vert &\leq &\frac{\left\vert
z\right\vert }{a}\left( 1+a\sum_{k=2}^{\infty }\left\vert c_{k}\right\vert
\left\vert z\right\vert ^{k-1}\right) \\
&\leq &\frac{\left\vert z\right\vert }{a}\left( 1+\sum_{k=2}^{\infty }\left( 
\frac{6}{a^{2}}\right) ^{k-1}\left\vert z\right\vert ^{k-1}\right) \\
&\leq &\frac{\left\vert z\right\vert }{a}\frac{1}{1-\frac{6}{a^{2}}%
\left\vert z\right\vert }\leq \frac{2}{a}\left\vert z\right\vert ,
\end{eqnarray*}%
where in the second line we used inequality %
\eqref{eq:coef_psi_inv_inequality}. QED.

\begin{lemma}
\label{lemma_psi_X_inverse3}Let $X$ be unitary and $EX=a>0.$ If $\left\vert
z\right\vert \leq a^{2}/12,$ and $\alpha \geq 1-a,$ then 
\begin{equation*}
\left\vert \frac{\psi _{X}^{-1}\left( z\right) }{z/\left( a+z\right) }%
-1\right\vert \leq \frac{3342}{a^{2}}\alpha \left\vert z\right\vert .
\end{equation*}
\end{lemma}

\textbf{Proof:} First of all, by Lemma \ref{lemma_psi_X_inverse}%
\begin{equation*}
\left\vert \psi _{X}^{-1}\left( z\right) \right\vert \leq \frac{2}{a}%
\left\vert z\right\vert
\end{equation*}%
for $\left\vert z\right\vert \leq a^{2}/12.$

Now we use the functional equation for $\psi _{X}^{-1}\left( z\right) $:%
\begin{equation*}
\psi _{X}\left( \psi _{X}^{-1}\left( z\right) \right) =z.
\end{equation*}%
If $\left\vert z\right\vert \leq a^{2}/12,$ then $\left\vert \psi
_{X}^{-1}\left( z\right) \right\vert \leq 2\left\vert z\right\vert /a\leq
a/6<1/2$ and we can apply Lemma \ref{lemma_psi_approximation} to get:%
\begin{eqnarray*}
\left\vert z-\frac{a\psi _{X}^{-1}\left( z\right) }{1-\psi _{X}^{-1}\left(
z\right) }\right\vert &\leq &716\alpha \left\vert \psi _{X}^{-1}\left(
z\right) \right\vert ^{2} \\
&\leq &716\alpha \frac{4}{a^{2}}\left\vert z\right\vert ^{2}=\frac{%
2864\alpha }{a^{2}}\left\vert z\right\vert ^{2}.
\end{eqnarray*}

Next, we write this as 
\begin{eqnarray*}
\left\vert z-\left( a+z\right) \psi _{X}^{-1}\left( z\right) \right\vert
&\leq &\left\vert 1-\psi _{X}^{-1}\left( z\right) \right\vert \frac{%
2864\alpha }{a^{2}}\left\vert z\right\vert ^{2} \\
&\leq &\frac{7}{6}\frac{2864\alpha }{a^{2}}\left\vert z\right\vert ^{2}<%
\frac{3342\alpha }{a^{2}}\left\vert z\right\vert ^{2}.
\end{eqnarray*}%
(In the second inequality we used the fact that $\left\vert \psi
_{X}^{-1}\left( z\right) \right\vert \leq 1/6$ if $\left\vert z\right\vert
\leq a^{2}/12.$) It follows that%
\begin{equation*}
\left\vert \frac{\psi _{X}^{-1}\left( z\right) }{z/\left( a+z\right) }%
-1\right\vert \leq \frac{3342}{a^{2}}\alpha \left\vert z\right\vert .
\end{equation*}%
QED.

\begin{lemma}
\label{lemma_product2}Let $EX_{i}=a_{i}$ and assume that for each $i,$ it is
true that $a_{i}\geq a.$ Assume also that $\left\vert z\right\vert \leq
a^{2}/3342$ and let $\alpha_{i}=:1-a_{i}$. Then 
\begin{equation*}
\left\vert \prod_{i=1}^{n}\frac{1}{\psi_{i}^{-1}\left(z\right)}\right\vert
\leq\frac{\prod\nolimits _{i=1}^{n}a_{i}}{\left\vert z\right\vert ^{n}}%
\prod_{i=1}^{n}\frac{1}{1-c_{i}\left\vert z\right\vert }\left\vert
\prod_{i=1}^{n}\left(1+\frac{z}{a_{i}}\right)\right\vert ,
\end{equation*}
where $c_{i}=3342\alpha_{i}/a_{i}^{2}.$
\end{lemma}

\textbf{Proof:} From Lemma \ref{lemma_psi_X_inverse3} we infer that 
\begin{equation*}
\left\vert \psi_{i}^{-1}\left(z\right)\right\vert \geq\left|\frac{z}{a_{i}}%
\frac{1}{1+z/a_{i}}\right|\left(1-\frac{3342\alpha_{i}}{a_{i}^{2}}\left\vert
z\right\vert \right).
\end{equation*}
Multiplying these inequalities together and inverting both sides, we get the
desired result. QED.

\begin{lemma}
Under the assumptions of the previous lemma, the following inequality holds: 
\begin{equation}
\left\vert f\left(z\right)\right\vert \leq\left(\left\vert 1+z\right\vert
\prod\nolimits _{i=1}^{n}a_{i}\prod_{i=1}^{n}\frac{1}{1-c_{i}\left\vert
z\right\vert }\left\vert \prod_{i=1}^{n}\frac{1+z/a_{i}}{1+z}\right\vert
\right)^{k},  \label{formula_fx_product_estimate}
\end{equation}
where $c_{i}=3342\alpha_{i}/a_{i}^{2}.$
\end{lemma}

\textbf{Proof:} The claim of this lemma is a direct consequence of Lemma \ref%
{lemma_product2} and equality (\ref{formula_product}). QED.

We will estimate terms in the product on the right-hand side of (\ref%
{formula_fx_product_estimate}) one by one.

\begin{lemma}
Suppose that $\alpha_{i}=:1-a_{i}<1-a$ for each $i,$ and that 
\begin{equation*}
\left\vert z\right\vert \leq\frac{a^{2}}{6684}\min\left\{
1,\left(\sum_{i=1}^{n}\alpha_{i}\right)^{-1}\right\} .
\end{equation*}
Then 
\begin{equation*}
\left\vert \prod_{i=1}^{n}\frac{1+z/a_{i}}{1+z}\right\vert \leq e.
\end{equation*}
\end{lemma}

\textbf{Proof:} We write:%
\begin{equation*}
\left\vert \prod_{i=1}^{n}\frac{1+z/a_{i}}{1+z}\right\vert =\exp \left( 
\mathrm{Re}\sum_{i=1}^{n}\log \left( 1+\frac{\alpha _{i}}{a_{i}}\frac{z}{1+z}%
\right) \right) .
\end{equation*}%
Recall that $\mathrm{Re}\log \left( 1+u\right) \leq \left\vert u\right\vert $
if $\left\vert u\right\vert <1.$ Under our assumption about $\left\vert
z\right\vert ,$ it is true that 
\begin{equation*}
\left\vert \frac{\alpha _{i}}{a_{i}}\frac{z}{1+z}\right\vert <1.
\end{equation*}%
Therefore we can write:%
\begin{eqnarray*}
\left\vert \prod_{i=1}^{n}\frac{1+z/a_{i}}{1+z}\right\vert &\leq &\exp
\left( \left\vert \frac{z}{1+z}\right\vert \sum_{i=1}^{n}\frac{\alpha _{i}}{%
a_{i}}\right) \\
&\leq &\exp \left( \frac{2}{a}\left\vert z\right\vert \sum \alpha _{i}\right)
\\
&\leq &e.
\end{eqnarray*}%
QED.

\begin{lemma}
Suppose that $\alpha_{i}=:1-a_{i}<1-a$ for each $i,$ and that 
\begin{equation*}
\left\vert z\right\vert \leq\frac{a^{2}}{6684}\min\left\{
1,\left(\sum_{i=1}^{n}\alpha_{i}\right)^{-1}\right\} .
\end{equation*}
Then,%
\begin{equation*}
\prod_{i=1}^{n}\frac{1}{1-c_{i}\left\vert z\right\vert }\leq e,
\end{equation*}
where $c_{i}=3342\alpha_{i}/a_{i}^{2}.$
\end{lemma}

\textbf{Proof:} We use the inequality $\log\left(1-u\right)\geq-2u,$ which
is valid for $u\in\left(0,1/2\right),$ and write:%
\begin{eqnarray*}
\prod_{i=1}^{n}\frac{1}{1-c_{i}\left\vert z\right\vert } & = &
\exp\left(-\sum_{i=1}^{n}\log\left(1-c_{i}\left\vert z\right\vert
\right)\right) \\
& \leq & \exp\left(2\left\vert z\right\vert \sum_{i=1}^{n}c_{i}\right) \\
& \leq & \exp\left[\frac{6684}{a^{2}}\left(\sum_{i=1}^{n}\alpha_{i}\right)%
\left\vert z\right\vert \right] \\
& \leq & e.
\end{eqnarray*}
QED.

Finally, note that if $\left\vert z\right\vert \leq a^{2}/6684,$ then $%
\left\vert 1+z\right\vert \leq2.$ Collecting all the pieces, we obtain that
if 
\begin{equation*}
\left\vert z\right\vert \leq\frac{a^{2}}{6684}\min\left\{
1,\left(\sum_{i=1}^{n}\alpha_{i}\right)^{-1}\right\}
\end{equation*}
then:%
\begin{equation*}
\left\vert f\left(z\right)\right\vert
\leq\left(2e^{2}\right)^{k}\left(\prod\nolimits _{i=1}^{n}a_{i}\right)^{k}.
\end{equation*}
This completes the proof of Proposition \ref{proposition_f_estimate}.

\section{Conclusion\label{section_conclusion}}

We have derived sufficient and necessary conditions for the product of free
unitary operators to converge in distribution to the uniform law. If \emph{%
essential convergence} denotes the situation when the partial products
continue to converge even after an arbitrary finite number of terms are
removed, then the necessary and sufficient condition for essential
convergence is that the products $\prod_{i=k_{0}}^{n}EX_{i}$ converges to
zero for all $k_{0},$ that is, that the products of expectations essentially
converge to zero. Essential convergence implies convergence. In addition,
non-essential convergence can occur when there is either a term that has the
uniform distribution, or there are two terms that have zero expectation. In
the latter case convergence occurs because the product of these two terms
has the uniform distribution.

\bibliographystyle{abbrvnat}
\bibliography{comtest}

\end{document}